\documentclass[preprint,12pt]{elsarticle}

\usepackage{amssymb}
\usepackage{xcolor}
\usepackage{mathrsfs}
\usepackage{amscd}
\usepackage{appendix}
\usepackage{geometry}
\usepackage{setspace}
\usepackage{amsfonts}
\usepackage{amsmath}
\usepackage{lineno,hyperref}
\usepackage{tikz}
\usetikzlibrary{matrix,arrows}

\usepackage{amssymb,amscd,amsthm,verbatim}

\biboptions{numbers,sort&compress}

\usepackage{bbm}

\date{\today}

\newtheorem{theorem}{Theorem}[section]
\newtheorem{lemma}[theorem]{Lemma}
\newtheorem{prop}[theorem]{Proposition}

\theoremstyle{definition}

\newtheorem*{definition 1}{Definition 1}
\newtheorem*{definition 2}{Definition 2}
\newtheorem*{definition 3}{Definition 3}
\newtheorem*{definition 4}{Definition 4}
\newtheorem*{definition 5}{Definition 5}
\newtheorem{example}[theorem]{Example}
\newtheorem{remark}[theorem]{Remark}

\theoremstyle{plain}

\allowdisplaybreaks \numberwithin{equation}{section}

\def\Proof{
\noindent \it Proof.\ \ \rm}

\journal{arXiv}

\begin{document}
\title{Construction and Characterization of Oscillatory Chain Sequences}

\author[1]{Zejun Dai}
\ead{daizejun@aliyun.com}

\author[2]{Daxiong Piao}
\ead{dxpiao@ouc.edu.cn}
\cortext[cor1]{Corresponding author}

\author[2]{Jinglai Qiao \corref{cor1}}
\ead{oucqiaojinglai@sohu.com}

\address[1]{Department of Applied Mathematics,  Liaonig Petrochemical University, Funshun 113001, P.R.China}
\address[2]{School of Mathematical Sciences,  Ocean University of China,Qingdao 266100, P.R.China}
\begin{abstract}

This paper initiates a theoretical investigation of $\frac{1}{4}$-oscillatory chain sequences $\{a_n\}$, generalizing Szwarc's classical framework for non-oscillatory chains \cite{Sz94, Sz98, Sz02, Sz03} to sequences fluctuating around $\frac{1}{4}$. We prove the existence of a fixed point for the critical map $f(x)=1-\frac{1}{4x}$ and establish convergence properties linking oscillatory behavior to parameter sequences $\{g_n\}$. A complete characterization is provided via a necessary and sufficient condition, exemplified by explicit solutions $a_n=\frac{1}{4}\left(1+(-1)^{n}\varepsilon_{n}\right)$. Crucially, we construct oscillatory chain sequences for which the series $\sum_{n=1}^{\infty} \left(a_n - \frac{1}{4}\right)$ diverges, demonstrating fundamentally different behavior outside the hypothesis $a_n \ge \frac{1}{4}$ required by Chihara's bound.


\end{abstract}

\maketitle


\section{Introduction}
\label{sec:intro}
Chain sequences constitute a fundamental analytical tool in both orthogonal polynomial theory and continued fraction analysis, as evidenced by seminal monographs \cite{Chi-book, Wall-book}.
Formally defined, a sequence $\{a_n\}_{n=1}^{\infty}$ qualifies as a chain sequence if it admits the decomposition:
\begin{equation}
a_n=g_n(1-g_{n-1}),\ \ \ n\geq 1,
\end{equation}
where the parameter sequence $\{g_n\}_{n=0}^{\infty}$ satisfies $0\leq g_n\leq 1$. It is worth noting that while we adopt the parameter sequence restrictions $0\leq g_n\leq 1$ following Wall \cite{Wall-book}, Chihara \cite[p. 91]{Chi-book} utilizes a more restrictive definition, namely $0\leq g_0 < 1$ and $0 < g_n < 1$ for $n\geq 1$. Since we consider a conjecture from Chihara's book, we remark this distinction, though it does not affect the validity of our main results.

The critical threshold of $\frac{1}{4}$ emerges as pivotal in chain sequence analysis; see \cite[Chapter III]{Chi-book}. While constant sequences $\{a\}$ are chain sequences if and only if $0 < a  \leq \frac{1}{4}$, Chihara's landmark conjecture proposed that for chain sequences with $a_n \geq \frac{1}{4}$, the series sum satisfies:
\begin{equation}\label{conjecture}
\sum_{n=1}^{+\infty} (a_n-\frac{1}{4})\leq \frac{1}{4}.
\end{equation}

Szwarc's seminal 1998 proof \cite{Sz98} not only confirmed this conjecture but also established refined estimates for parameter-dependent sequences. Subsequent work \cite{Sz02} demonstrated the maximality of sequences like $a_n=\frac{1}{2}+\frac{1}{4(4n^2-1)}$, crucially informing Jacobi operator norm bounds, while \cite{Sz03} further established sharp parameter-dependent estimates through refined continued fraction techniques. These developments, however, remained confined to non-oscillatory sequences satisfying $a_n \geq \frac{1}{4}+b_n\left(b_n>0\right)$, leaving open the oscillatory regime where $a_n$ fluctuates about $\frac{1}{4}$.

This paper bridges this theoretical gap by extending Szwarc's framework to $\frac{1}{4}$-oscillatory chain sequences where both $P_{\frac{1}{4}}=\left\{k: a_k>\frac{1}{4}\right\}$ and $N_{\frac{1}{4}}= \left\{k: a_k<\frac{1}{4}\right\}$ are infinite. Our key contributions include:
\begin{enumerate}[(1)]
\item \emph{Fixed Point Analysis}: Establishing existence and convergence properties for the critical mapping $f(x)=1-\frac{1}{4x}$, linking oscillatory dynamics to parameter sequences.

\item \emph{Characterization Theorem}: Deriving necessary and sufficient conditions for $\frac{1}{4}$-oscillatory chain sequences, exemplified by explicit constructions of the form $a_n=\frac{1}{4}\left(1+(-1)^n \varepsilon_n\right)$.

\item \emph{Divergence Phenomena}: Constructing divergent oscillatory chain sequences that exhibit divergent summation of $\sum_{n=1}^{\infty} \left(a_n - \frac{1}{4}\right)$. This demonstrates fundamentally different behavior outside the hypothesis $a_n \geq \frac{1}{4}$ required by Chihara's bound, revealing fundamental limitations of prior estimates when applied to oscillatory regimes.

\end{enumerate}

The paper proceeds as follows: Section \ref{sec:2} formalizes oscillatory chain sequences and establishes their convergence properties. Section \ref{sec:exam} presents concrete examples illustrating this divergent behavior outside the scope of Chihara's condition. Our results collectively extend the analytical toolkit of \cite{Sz94, Sz98, Sz02, Sz03} while uncovering novel phenomena unique to oscillatory regimes.

\section{Oscillatory chain sequences}\label{sec:2}

For a sequence $\{a_n\}_{n=1}^{\infty}$, we define two sets as
$$ P_a=\{k: a_k-a>0, k\in \mathbb{N}\}$$
 and
 $$ N_a=\{k: a_k-a<0, k\in \mathbb{N}\}.$$
 If both $P_a$ and $N_a$ are infinite, then we say that $\{a_n\}_{n=1}^{\infty}$ is an $a$-oscillatory sequence, or that $\{a_n\}_{n=1}^{\infty}$ oscillates around $a$. Otherwise, it is a non-$a$-oscillatory sequence. If $\{a_n\}_{n=1}^{\infty}$ does not oscillate around any number, then we simply say that $\{a_n\}_{n=1}^{\infty}$ is non-oscillatory. In this sense, one could say that Szwarc \cite{Sz94, Sz98, Sz02, Sz03} discussed mainly the non-oscillatory chain sequences, but not oscillatory chain sequences.

It is easy to see that if $\{a_n\}_{n=1}^{\infty}$ is a chain sequence with a parameter sequence $\{g_n\}_{n=1}^{\infty}$, then the existence of $\lim_{n\rightarrow\infty} g_n$ implies that $\lim_{n\rightarrow\infty} a_n$ exist, while the converse is generally not true. But we can prove following proposition.

\begin{prop}\label{prop2-1}
 Suppose $\lim_{n\rightarrow \infty}a_n = \frac{1}{4}$. If $\{a_n\}$ is a chain sequence with a parameter sequence $\{g_n\}$, then $\lim_{n\rightarrow \infty}g_n$ exists and  $\lim_{n\rightarrow \infty}g_n=\frac{1}{2}$.\\
\end{prop}

To prove this proposition, we need following lemma.

\begin{lemma}\label{lem2-1}
Suppose $f(x)=1-\frac{1}{4x}, \ \ x\in[\frac{1}{2},1]$. Define a sequence $x_0, x_1,\cdots, x_n,\cdots $ by iteration $x_{i+1}=f(x_i),i=0,1,2 \dots$,
then $\lim_{n\rightarrow \infty} x_n= \frac{1}{2}$.

\proof  It is obvious that $f(x)$ is monotonically increasing on $[\frac{1}{2}, 1]$ and $f(x)\in[\frac{1}{2},1]$ for $x\in [\frac{1}{2}, 1]$. Because of $f(x)-x=-\frac{(2x-1)^2}{4x}\le 0$, sequence $x_{i+1}=f(x_i), i=0,1,2,\dots$ is non-increasing and has a lower bound of $\frac{1}{2}$. Hence $\displaystyle\lim_{i\to \infty}x_{i}$ exists.\\

Taking the limit on both sides of the equation $x_{i+1}=1-\frac{1}{4x_i}$, we obtain $$\displaystyle\lim_{i\to \infty}x_{i}=\frac{1}{2}.$$
\qed
\end{lemma}

{\it Proof of Proposition \ref{prop2-1}}\ \   Since $a_n\rightarrow\frac{1}{4}$, for any $\epsilon>0$, there exists an integer  $N>0$ , such that the inequality
\begin{equation*}
g_{n+k}(1-g_{n+k-1})\ge\frac{1-\epsilon}{4}
\end{equation*}
holds for any $n>N$ and $ k>\mathbb{N}_0$. As a result,
\begin{equation*}
g_{n+k-1}\le 1-\frac{1-\epsilon}{4g_{n+k}}.
\end{equation*}
Based on the arbitrariness of $\epsilon$, we have
\begin{equation}
g_{n+k-1}\le 1-\frac{1}{4g_{n+k}}.\label{1}
\end{equation}
This implies that for $n>N$,
\begin{equation}
g_n\le f(g_{n+1})\le f\circ f(g_{n+2})\le \dots\le f^k(g_{n+k})\le f^k(1),
\end{equation}
where $f(x)=1-\frac{1}{4x}$, and $f^k$ denotes the $k$-th composition of function $f$.

Let $k\to \infty$, then Lemma \ref{lem2-1} yields $f^k(1) \to \frac{1}{2}$, and so $g_n\le \frac{1}{2}$ for $n\ge N$.
\par We assume that $g_n=\frac{1}{2}-\delta_n$, then $\delta_n\ge 0$ for $n\ge N$.
\par If $\delta_n> 0$ for $n\ge N$, then
\begin{equation}
a_{n+1}=g_{n+1}(1-g_n)=(\frac{1}{2}-\delta_{n+1})(1-(\frac{1}{2}-\delta_n))\ge \frac{1-\epsilon}{4},
\end{equation}
i.e.
\begin{equation}
(\frac{1}{2}-\delta_{n+1})(\frac{1}{2}+\delta_n)\ge \frac{1-\epsilon}{4}.\label{2.4}
\end{equation}
Further simplification yields
$\delta_n - \delta_{n+1} \geq -\frac{\epsilon}{2}.$

Assume $\limsup_{n \to \infty} \delta_n = \beta > 0$. Then there exists a subsequence $\delta_{n_k} \to \beta$. From the asymptotic relation derived from (\ref{2.4}):
$\delta_{n+1} = \frac{\delta_n}{1 + 2\delta_n} + o(1), \quad \text{as } n \to \infty,$
we iteratively obtain:
\begin{align*}
\delta_{n_k+1} &\to \frac{\beta}{1 + 2\beta}, \\
\delta_{n_k+2} &\to \frac{\beta}{1 + 4\beta}, \\
\delta_{n_k+m} &\to \frac{\beta}{1 + 2m\beta} \quad \text{for any fixed } m.
\end{align*}
For sufficiently large $m$, $\frac{\beta}{1 + 2m\beta} < \frac{\beta}{2}$. This contradicts $\beta$ being the limit superior. Hence $\beta = 0$, and $\delta_n \to 0$, we obtain that $\displaystyle\lim_{n\to \infty}g_{n}=\frac{1}{2}$.
\par If there exists $n_0>N$, such that $\delta_{n_0}=0$, then
\begin{equation}
a_{n_0+1}=(\frac{1}{2}-\delta_{n_0+1})(1-(\frac{1}{2}-0))\ge\frac{1-\epsilon}{4}.
\end{equation}
Thus $\delta_{n_0+1}\le \frac{\epsilon}{2}$, i.e. $\delta_{n_0+1}\to 0$, $g_{n_0+1}\to \frac{1}{2}.$ By induction, we still obtain that $\displaystyle\lim_{n\to \infty}g_{n}=\frac{1}{2}$.
\qed

\begin{theorem}\label{th2-3} Let $a_n=\frac{1}{4}(1+(-1)^n\varepsilon_n)$, with $\varepsilon_n\ge 0$ and  $\displaystyle\lim_{n\to \infty}\varepsilon_{n}=0$. Then $\{a_n\}$ is a chain sequence if and only if there exists a sequence $\{c_n\}$ of positive numbers such that
\begin{itemize}
\item[(i)]$c_{n+1}\le 2c_n;$
\item[(ii)]$ \sum\limits_{m=n}^{\infty}(-1)^{m-n}c_m\varepsilon_m=(-1)^n(c_{n+1}-c_n).$
\end{itemize}
\end{theorem}
\Proof (1). \emph{Necessity}.

Let $a_n=g_n(1-g_{n-1})$, $g_n=\frac{1}{2}(1+(-1)^n\delta_n)$. Then
\begin{equation*}
\frac{1}{4}(1+(-1)^n\varepsilon_n)=\frac{1}{2}(1+(-1)^n\delta_n)(1-\frac{1}{2}(1+(-1)^{n-1}\delta_{n-1})).
\end{equation*}
Moreover,
\begin{equation}
\varepsilon_n=\delta_{n-1}+\delta_n+(-1)^n\delta_{n-1}\delta_n.\label{epsilon_n}
\end{equation}
Set $c_1=1$ and
\begin{equation}
\frac{c_{n+1}}{c_n}=1+(-1)^n\delta_{n-1}.\label{c_n}
\end{equation}
Then
\begin{equation*}
c_{n+1}-c_n=(-1)^n\delta_{n-1}c_n
\end{equation*}
and $c_{n+1}\le 2c_n$, since $\delta_n\le 1$. According to (\ref{epsilon_n}),
\begin{equation}
\varepsilon_n=\delta_{n-1}+\delta_n[1+(-1)^n\delta_{n-1}]=\delta_{n-1}+\delta_n\frac{c_{n+1}}{c_n}.
\end{equation}
Multiplying both sides of the equation above by $c_n$ gives
\begin{equation}
c_n\varepsilon_n=c_n\delta_{n-1}+\delta_n c_{n+1}.
\end{equation}
Hence,
\begin{equation}
(-1)^{m-n} c_m\varepsilon_m=(-1)^{m-n} c_m\delta_{m-1}+(-1)^{m-n} c_{m+1}\delta_m.\label{-1c_nepsilon_n}
\end{equation}
Summing up (\ref{-1c_nepsilon_n}) yields that
\begin{equation}
\sum\limits_{m=n}^{\infty}(-1)^{m-n} c_m\varepsilon_m=c_n\delta_{n-1}+\displaystyle\lim_{m\to \infty}(-1)^{m-n} c_{m+1}\delta_m.
\end{equation}
According to Proposition \ref{prop2-1}, $\displaystyle\lim_{m\to \infty}\delta_m=0$. We can further see that $\{ c_{n}\}$ is bounded by (\ref{c_n}). As a result $\displaystyle\lim_{m\to \infty}(-1)^{m-n} c_{m+1}\delta_m=0$. Then
\begin{equation}
\sum\limits_{m=n}^{\infty}(-1)^{m-n} c_m\varepsilon_m=c_n\delta_{n-1}=(-1)^n(c_{n+1}-c_n).
\end{equation}

(2). \emph{Sufficiency.}

Set $$\delta_{n-1}=\frac{1}{c_n}\sum\limits_{m=n}^{\infty}(-1)^{m-n} c_m\varepsilon_m,$$
 then
\begin{equation*}
|\delta_n|=\frac{1}{c_{n+1}}|\sum\limits_{m=n+1}^{\infty}(-1)^{m-n} c_m\varepsilon_m|=\frac{1}{c_{n+1}}|c_{n+1}-c_n|=|1-\frac{c_n}{c_{n+1}}|.
\end{equation*}
By (i), we have $|1-\frac{c_n}{c_{n+1}}|\le \frac{1}{2}$. Hence $|\delta_n|\le \frac{1}{2}$. By (ii), we have
\begin{equation*}
\varepsilon_n=\delta_{n-1}+\frac{c_{n+1}}{c_n}\delta_n,
\end{equation*}
thus
\begin{equation}
(-1)^n\delta_{n-1}=(-1)^n\varepsilon_n-(-1)^n\frac{c_{n+1}}{c_n}\delta_n\label{-1delta_n}
\end{equation}
\par Let $h_n=\frac{1}{2}(1+(-1)^n\delta_n)$, then
\begin{equation}
4h_n(1-h_{n-1})=1+(-1)^n\delta_{n-1}+(-1)^n\delta_n+\delta_n\delta_{n-1}.\label{4bn}
\end{equation}
Substitute (\ref{-1delta_n}) into (\ref{4bn}), we can obtain
\begin{equation}
\begin{split}
4h_n(1-h_{n-1})&=1+(-1)^n\varepsilon_n-(-1)^n\frac{c_{n+1}}{c_n}\delta_n+(-1)^n\delta_n+\delta_n\delta_{n-1}\\
&=1+(-1)^n\varepsilon_n+(-1)^n\delta_n\frac{c_n-c_{n+1}}{c_n}+\delta_n\delta_{n-1}.\\\label{4bn2}
\end{split}
\end{equation}
Substitute (ii) into (\ref{4bn2}), we have
\begin{equation*}
\begin{split}
4h_n(1-h_{n-1})&=1+(-1)^n\varepsilon_n+(-1)^n\delta_n(-1)^{n+1}\frac{1}{c_n}(\sum\limits_{m=n}^{\infty}(-1)^{m-n}c_m\varepsilon_m)+\delta_n\delta_{n-1}\\
&=1+(-1)^n\varepsilon_n-\frac{1}{c_n}(\sum\limits_{m=n}^{\infty}(-1)^{m-n}c_m\varepsilon_m)\delta_n+\delta_n\delta_{n-1}\\
&=1+(-1)^n\varepsilon_n-\delta_{n-1}\delta_n+\delta_n\delta_{n-1}\\
&=1+(-1)^n\varepsilon_n\\
&=4a_n
\end{split}
\end{equation*}
Thus, $\{a_n\}$ is a chain sequence.
\qed

\begin{remark}\label{rem:con-exam}
In Theorem \ref{th2-3}, we suppose that $\varepsilon_n\ge0$. If we require that $\varepsilon_n>0$, then the sequence $\{a_n=\frac{1}{4}(1+(-1)^n\varepsilon_n)\}$ is a $\frac{1}{4}-$ oscillatory one.
\end{remark}

\begin{theorem}\label{thm:p-q}
Let
\begin{equation}\label{con:p-q}
a_{2k-1}=\frac{1}{4}-p,\quad  a_{2k}=\frac{1}{4}+q, \quad  k=1,2,\cdots, \quad p,q>0.
\end{equation}
 If there exist $\epsilon,\gamma \in (0,1)$, $\epsilon<\gamma$, such that
 \begin{equation}\label{con:eps-gam}
 \epsilon \leq \frac{1}{4}-p\leq \gamma(1-\gamma),\,\, \epsilon \leq \frac{1}{4}+q\leq \gamma(1-\frac{\frac{1}{4}-p}{1-\gamma})= \gamma-\frac{\gamma}{1-\gamma}(\frac{1}{4}-p),
 \end{equation}
 then $\{a_{n}\}$ is a chain sequence and $\sum_{n=1}^{+\infty} (a_n-\frac{1}{4})$ is divergent.
\end{theorem}
\proof We prove the existence of the corresponding parameter sequence by construction.

Fix $g_0=0$, then $g_1=a_1=\frac{1}{4}-p\in[\epsilon, \gamma(1-\gamma)]\subset[\epsilon,\gamma]\subset[0, 1]$, and $$g_2=\frac{a_2}{1-g_1}=\frac{\frac{1}{4}+q}{1-(\frac{1}{4}-p)}\in[\epsilon, \gamma]\subset[0, 1].$$ By this way, we obtain a sequence $ \{g_n\}$.

In the following, we verify that $0 \leq g_n \leq 1$ .

Suppose that $g_{2k-1}, g_{2k}\in [\epsilon, \gamma]$ when $k\leq n$, then we have that
\begin{equation*}
\begin{split}
&g_{2k+1}=\frac{a_{2k+1}}{1-g_{2k}}=\frac{\frac{1}{4}-p}{1-g_{2k}}\geq \frac{\frac{1}{4}-p}{1-\epsilon}\geq\frac{\epsilon(1-\epsilon)}{1-\epsilon}= \epsilon,\\
&g_{2k+1}=\frac{a_{2k+1}}{1-g_{2k}}=\frac{\frac{1}{4}-p}{1-g_{2k}}\leq \frac{\frac{1}{4}-p}{1-\gamma}\leq\frac{\gamma(1-\gamma)}{1-\gamma}= \gamma,\\
&g_{2k+2}=\frac{a_{2k+2}}{1-g_{2k+1}}=\frac{\frac{1}{4}+q}{1-g_{2k+1}}\geq \frac{\frac{1}{4}+q}{1-\epsilon}\geq\frac{\epsilon(1-\epsilon)}{1-\epsilon}= \epsilon,\\
&g_{2k+2}=\frac{a_{2k+2}}{1-g_{2k+1}}=\frac{\frac{1}{4}+q}{1-g_{2k+1}}=\frac{\frac{1}{4}+q}{1-\frac{\frac{1}{4}-p}{1-g_{2k}}}\leq \frac{\gamma(1-\frac{\frac{1}{4}-p}{1-\gamma})}{1-\frac{\frac{1}{4}-p}{1-\gamma}}=\gamma.\\
\end{split}
\end{equation*}
And so we conclude that $g_{2k+1}, g_{2k+2}\in [\epsilon,\gamma]$. By induction, $g_n\in  [\epsilon,\gamma]\subset (0,1)$. As a result, $\{a_n\}$ is a chain sequence. It is obvious that $\sum_{n=1}^{+\infty} (a_n-\frac{1}{4})=-p+q-p+q-\cdots$ is divergent.
\qed
\remark Because $\frac{1}{4}-p\leq \gamma(1-\gamma)$, it holds that $\gamma-\frac{\gamma}{1-\gamma}(\frac{1}{4}-p)\geq\gamma(1-\gamma)$. Taking $\gamma-\frac{\gamma}{1-\gamma}(\frac{1}{4}-p)$ as the upper bound of $a_{2k}$ instead of $\gamma(1-\gamma)$ can avoid the situation where $a_{2k}$ and $g_{2k}$ belongs to empty set.

\section{Examples}\label{sec:exam}

While the chain sequences themselves may converge or diverge, the key feature outside the scope of Chihara's bound stems from the divergence of $\sum (a_n - \frac{1}{4})$.

\subsection{Convergent Oscillatory Chain Sequence}  
\begin{example}

\par A sequence ${\{a_n\}}$ is given as
\begin{equation}\label{eq:os-sequence1}
a_n=\frac{1}{4}+(-1)^n\frac{1}{4(4n^2-1)}.
\end{equation}
It is obvious that ${\{a_n\}}$ oscillates around $\frac{1}{4}$.
Let \begin{equation*}
g_n=\frac{2n^2-\frac{1}{4}+(-1)^n\frac{1}{4}}{(2n+1)^2},
\end{equation*}
then we have
\begin{equation*}
\begin{split}
g_n(1-g_{n-1})&=\frac{2n^2-\frac{1}{4}+(-1)^n\frac{1}{4}}{(2n+1)^2}\left[1-\frac{2(n-1)^2-\frac{1}{4}+(-1)^{(n-1)}\frac{1}{4}}{(2n-1)^2}\right]\\
&=\frac{2n^2-\frac{1}{4}+(-1)^n\frac{1}{4}}{(2n+1)^2}\frac{(2n-1)^2-2(n-1)^2+\frac{1}{4}+(-1)^n\frac{1}{4}}{(2n-1)^2}\\
&=\frac{1}{(4n^2-1)^2}\left[2n^2-\frac{1}{4}+(-1)^n\frac{1}{4}\right]\left[2n^2-\frac{3}{4}+(-1)^n\frac{1}{4}\right]\\
&=\frac{1}{(4n^2-1)^2}\left[4n^4-2n^2+\frac{1}{4}+(-1)^n(n^2-\frac{1}{4})\right]\\
&=\frac{1}{4}+(-1)^n\frac{1}{4(4n^2-1)}=a_n.
\end{split}
\end{equation*}
It is easy to verify that both ${g_n}$ and $a_n$ belong to $[0,1]$ for any positive integer $n$,  so ${\{a_n\}}$ is a $\frac{1}{4}$-oscillatory chain sequence and $\{g_n\}$ is the corresponding parameter sequence.
\par On the one hand, $a_n$ tends to $\frac{1}{4}$, as a result  $g_n$ tends to $\frac{1}{2}$. On the other hand, let
$$c_m=2^{m-1}\cdot c_1\cdot\prod_{i=1}^{m-1}\frac{[1+(-1)^i](2i^2+1)^2}{[8i^2-1+(-1)^i](4i^2-1)},$$
then $c_{m+1}\leq2c_m$ and$$\sum_{m=1}^{\infty}(-1)^{m-1}c_m\varepsilon_m=-1(c_2-c_1).$$
This verifies Proposition \ref{prop2-1} and Theorem \ref{th2-3}.

Additionally, in this case, $\sum_{n=1}^{+\infty} (a_n-\frac{1}{4})$ is convergent and the inequality (\ref{conjecture}) holds, since
\begin{align*}
\sum_{n=1}^{+\infty} |a_n-\frac{1}{4}|
&=\lim_{k\to+\infty}\sum_{n=1}^{k} |a_n-\frac{1}{4}|\\
&=\lim_{k\to+\infty}\sum_{n=1}^{k} |(-1)^n\frac{1}{4(4n^2-1)}|\\
&=\lim_{k\to+\infty}\sum_{n=1}^{k}\frac{1}{8}\left(\frac{1}{2n-1}-\frac{1}{2n+1}\right)
\end{align*}
\begin{align*}
&=\lim_{k\to+\infty}\frac{1}{8}\left(1-\frac{1}{2k+1}\right)\\
&=\frac{1}{8}< \frac{1}{4}.
\end{align*}

\end{example}

\begin{remark}
These results establish a dichotomy: while some oscillatory chain sequences satisfy the inequality of Chihara's bound \eqref{conjecture} (e.g., Example \ref{eq:os-sequence1}), others exhibit divergent summation since they do not satisfy the condition $a_n \geq \frac{1}{4}$ required by Chihara's hypothesis (e.g., Example \ref{exam:p-q}). This behavioral divergence distinguishes them from non-oscillatory chains above $\frac{1}{4}$, which uniformly satisfy \eqref{conjecture} \cite{Sz98}. Crucially, oscillatory chains preserve convergence properties (Proposition \ref{prop2-1}) while exhibiting novel phenomena absent in their non-oscillatory counterparts.
\end{remark}

\subsection{Divergent Oscillatory Chain Sequence}  

\begin{example}\label{exam:p-q}
 We present a concrete realization of the chain sequence in Theorem \ref{thm:p-q}. Set $p = 0.22$ and $q = 0.24$, giving:
$$ a_{2k-1} = 0.03, \quad a_{2k} = 0.49. $$

Select $\epsilon = 0.02$ and $\gamma = 0.6$. Then:
\begin{align*}
\frac{1}{4} - p &= 0.25 - 0.22 = 0.03, \\
\frac{1}{4} + q &= 0.25 + 0.24 = 0.49, \\
\gamma(1-\gamma) &= 0.6 \times 0.4 = 0.24, \\
\gamma - \frac{\gamma}{1-\gamma}\left(\frac{1}{4}-p\right) &= 0.6 - \frac{0.6}{0.4} \times 0.03 = 0.6 - 1.5 \times 0.03 = 0.555.
\end{align*}
Thus we have:
\begin{align*}
\epsilon = 0.02 &< 0.03 = \frac{1}{4} - p \leq 0.24 = \gamma(1-\gamma), \\
\epsilon = 0.02 &< 0.49 = \frac{1}{4} + q \leq 0.555 = \gamma - \frac{\gamma}{1-\gamma}\left(\frac{1}{4}-p\right)
\end{align*}
confirming condition \eqref{con:eps-gam} holds.

 Fix $g_0 = 0$, then iteratively define:
$$ g_1 = a_1 = 0.03, \quad g_2 = \frac{a_2}{1-g_1} = \frac{0.49}{1-0.03} = \frac{49}{97} \approx 0.505. $$

For $k \geq 1$, assume $g_{2k-1}, g_{2k} \in [0.02, 0.6]$, then
\begin{equation*}
\begin{split}
&g_{2k+1}=\frac{a_{2k+1}}{1-g_{2k}}=\frac{0.03}{1-g_{2k}}\geq \frac{0.03}{1-0.02}=\frac{3}{98}\geq 0.02,\\
&g_{2k+1}=\frac{a_{2k+1}}{1-g_{2k}}=\frac{0.03}{1-g_{2k}}\leq \frac{0.03}{1-0.6}=\frac{3}{40}\leq 0.6,\\
&g_{2k+2}=\frac{a_{2k+2}}{1-g_{2k+1}}=\frac{0.49}{1-g_{2k+1}}\geq \frac{0.49}{1-0.02}=\frac{49}{98}\geq 0.02,\\
&g_{2k+2}=\frac{a_{2k+2}}{1-g_{2k+1}}=\frac{0.49}{1-g_{2k+1}}=\frac{0.49}{1-\frac{0.03}{1-g_{2k}}}\leq \frac{0.49}{1-\frac{0.03}{1-0.6}}=\frac{49}{925}\leq 0.6,\\
\end{split}
\end{equation*}
that is $g_{2k+1}, g_{2k+2}\in [0.02,0.6]$.  By induction, we obtain that $g_n\in [0.02,0.6]\subset (0,1)$ for all $n$, confirming $\{a_n\}$ is a chain sequence.

We also have that $a_{2k-1}-\frac{1}{4}=-0.22$, $a_{2k}-\frac{1}{4}=0.24$. So the alternating series $\sum_{n=1}^{+\infty} (a_n-\frac{1}{4})$ is divergent. This shows that the inequality in (\ref{conjecture}) does not hold for these oscillatory chain sequences, which highlights their distinct behavior outside Chihara's condition.
\end{example}

\begin{remark}
Example \ref{exam:p-q} yields that the condition \eqref{con:eps-gam} is not empty.
\end{remark}

\vskip1cm

\section*{Acknowledgments}
The authors would like to thank the anonymous reviewer for their careful reading of the manuscript and constructive comments.

We thank Bo Wang and Zehong Wang for valuable discussions.

D.P and J.Q were supported in part by the  NSFC (No. 11571327, 11971059).


\vskip5mm
\section*{References}

\vskip5mm
\begin{thebibliography}{00}


\bibitem{Chi-book} T. S. Chihara, An Introduction to Orthogonal Polynomials, Gordon and Breach, New York, 1978.
\bibitem{Chi-trans} T. S. Chihara, Chain sequences and orthogonal polynomials, \emph{Trans. Amer. Math. Soc.} 104.
(1962), 1-16.
\bibitem{Sz94} R. Szwarc,  Chain sequences and compact perturbations of orthogonal polynomials.\emph{ Math. Z.} \textbf{217(}1994), 57 - 71.

\bibitem{Sz98} R. Szwarc,  Chain sequences, orthogonal polynomials, and Jacobi matrices, \emph{J. Approx. Theory,}  \textbf{92}(1998), 59 - 73.

\bibitem{Sz02} R. Szwarc,  Sharp estimates for Jacobi matrices and chain sequences,\emph{ J. Approx. Theory,} \textbf{118}(2002), 94 - 105.

\bibitem{Sz03} R. Szwarc,  Sharp estimates for Jacobi matrices and chain sequences, II.\emph{ J. Approx. Theory,} \textbf{125} (2003), 295-302.

\bibitem{Wall-book} H. S.  Wall,  Analytic Theory of Continued Fractions,  D. van Nostrand Co., New York, 1948.

\end{thebibliography}
\end{document}